

\input amstex
\documentstyle{amsppt}
\refstyle{A}
\magnification=\magstephalf
\nologo
\hsize=5.9in



\topmatter
\title
Loop coproducts in string topology and triviality of higher genus TQFT operations
\endtitle
\rightheadtext{Loop coproducts and higher genus TQFT string operations}
\author
Hirotaka Tamanoi
\endauthor
\affil
University of California, Santa Cruz
\endaffil
\address
Department of Mathematics, University of California Santa Cruz, 
Santa Cruz, CA 95064
\endaddress
\email
tamanoi\@math.ucsc.edu
\endemail
\keywords
loop coproduct, loop product, string operations, string topology
\endkeywords
\subjclass
55P35
\endsubjclass
\abstract
Cohen and Godin constructed positive boundary topological quantum
field theory (TQFT) structure on the homology of free loop spaces of oriented closed smooth manifolds by associating a certain operations called string operations to orientable surfaces with parametrized boundaries. We show that all TQFT string operations associated to surfaces of genus at least one vanish identically. This is a simple consequence of properties of the loop coproduct which will be discussed in detail. One interesting property is that  the loop coproduct is nontrivial only on the degree $d$ homology group of the connected component of $LM$ consisting of contractible loops, where $d=\dim M$, with values in the degree $0$ homology group of constant loops. Thus the loop coproduct behaves in a dramatically simpler way than the loop product. 
\endabstract
\toc
\head
1. Introduction and triviality of higher genus TQFT string operations  
\dotfill 1
\endhead
\head
2. The loop coproduct and its Frobenius compatibility \dotfill 4
\endhead
\head
3. Properties of the loop coproduct and their consequences \dotfill 11
\endhead
\head
{} References \dotfill 15
\endhead
\endtoc
\endtopmatter

\document

\specialhead
\S 1.\quad Introduction and triviality of higher genus TQFT
string operations
\endspecialhead

Let $M$ be a connected closed orientable smooth manifold of dimension
$d$, and let $LM=\text{Map}(S^1,M)$ be its free loop space of
continuous maps from the circle $S^1$ to $M$. Chas and Sullivan
\cite{CS} showed that its homology $\Bbb H_*(LM)=H_{*+d}(LM)$ comes
equipped with an associative graded commutative product of degree
$-d$, and a compatible Lie bracket of degree $1$.  These two products
together with an operator $\Delta$ of degree $1$ with $\Delta^2=0$,
coming from the natural $S^1$ action on $LM$, give $\Bbb H_*(LM)$ the
structure of a Batalin-Vilkovisky algebra.

The associative product called the loop product was generalized to so
called string operations by Cohen and Godin \cite{CG}. Let $\Sigma$ be an orientable
connected surface of genus $g$ with $p$ incoming and $q$ outgoing
parametrized boundary circles, where we require that $q\ge1$. To such
a surface $\Sigma$, they associated an operator $\mu_{\Sigma}$ of the
form
$$
\mu_{\Sigma}: H_*\bigl((LM)^p\bigr) @>>>
H_{*+\chi(\Sigma)d}\bigl((LM)^q\bigr),
$$ 
in such a way that $\mu_{\Sigma}$ depends only on the topological type of
the surface $\Sigma$ and $\mu_{\Sigma}$ is compatible with sewing of surfaces
along parametrized boundaries. These operations give rise to topological quantum field theory (TQFT) without a counit. When $\Sigma$ is a pair of pants with
either 2 incoming or 2 outgoing circles, we get a product and
and a coproduct:
$$
\aligned
\mu &: H_*(LM\times LM) \longrightarrow H_{*-d}(LM),\\
\Psi &: H_*(LM) \longrightarrow H_{*-d}(LM\times LM),
\endaligned
$$
where the product $\mu$ coincides with the loop product of Chas and
Sullivan. See formula (2-8) for a homotopy theoretic definition of the
loop product. Since any surface $\Sigma$ can be decomposed into pairs
of pants and capping discs, we can compute the string operation
$\mu_{\Sigma}$ by composing loop products and loop
coproducts according to pants decompositions of $\Sigma$. In
this paper, we study properties of coproduct in detail, and as a consequence we show that for higher genus surfaces $\Sigma$, the
string operations $\mu_{\Sigma}$ are always trivial.

\proclaim{Theorem A} Let $\Sigma$ be an oriented connected compact 
surface of genus $g$ with $p$ incoming and $q\ge1$ outgoing
parametrized boundary circles.  If $g\ge1$ or $q\ge3$, then the
associated string operation $\mu_{\Sigma}$ vanishes.
\endproclaim 

Thus the only nontrivial TQFT string operations correspond to genus $0$
surfaces with at most $2$ outgoing circles. To elements $a_1, a_2,
\dots, a_p\in H_*(LM)$, such operations associate either their loop 
product $a_1a_2\cdots a_p$ or its loop coproduct $\Psi(a_1a_2\cdots
a_p)$. Thus once we understand the loop coproduct $\Psi$, we know the
behavior of all string operations associated to orientable surfaces
with parametrized boundaries.  For $a\in H_*(LM)$, let $|a|$ denote
its homological degree. Let $c_0$ be the constant loop at the base
point $x_0$ in $M$, and let $[c_0]$ its homology class in $H_0(LM)$. 

The connected components of $LM$ are parametrized by conjugacy classes of $\pi_1(M)$. Let $(LM)_{[1]}$ be the component corresponding to the conjugacy class of $1\in\pi_1(M)$. This is the space of contractible loops in $M$. 

In addition to the Frobenius compatibility (Theorem 2.2), 
properties of the loop coproduct are described in Theorem B, whose part (2) shows dramatic simplicity of the loop coproduct compared with the loop product. Theorem B is the main result of this paper. Theorem A is only one of the consequences of Theorem B. We will discuss some of the other consequences in Theorem C.

\proclaim{Theorem B} \rom{(1)} Let $p\ge0$, and let $a_1,a_2,\dots, a_p\in
H_*(LM)$ be arbitrary $p$ elements. Then the image of the loop
coproduct $\Psi$ lies in the subset $H_*(LM)\otimes H_*(LM)\subset
H_*(LM\times LM)$ of cross products, and for any $0\le \ell\le p$ it
is given by
$$
\Psi(a_1\cdot a_2\cdots a_p)=\chi(M)[c_0]a_1\cdot a_2\cdot\cdots\cdot
a_{\ell}
\otimes[c_0]a_{\ell+1}\cdot\cdots\cdot a_p\in H_*(LM)\otimes H_*(LM),
$$
where $\chi(M)$ is the Euler characteristic of $M$. 

\rom{(2)} The loop coproduct $\Psi$ is nontrivial only on 
$H_d\bigl((LM)_{[1]}\bigr)$, the degree $d$ homology group of the component of contractible loops in $M$. On $H_d\bigl((LM)_{[1]}\bigr)$, the loop coproduct $\Psi$ has values in the homology classes of constant loops 
$H_0\bigl((LM)_{[1]}\bigr)\otimes H_0\bigl((LM)_{[1]}\bigr)
\cong\Bbb Z[c_0]\otimes[c_0]$. 
\endproclaim

Theorem B is proved in Theorem 3.1. 
Note that if $M$ has vanishing Euler characteristic, for
example if $M$ is odd dimensional, then its loop coproduct is
identically $0$. Before we prove the above result in \S3, in \S2 we will prove various general results on the loop 
coproduct including Frobenius compatibility (Theorem 2.2) with precise treatment of signs, and Frobenius compatibility and coderivation compatibility with respect to cap products (Theorem 2.4). 

Since the proof of Theorem A is more or less straightforward, we give its proof below. This vanishing property is the basis of triviality of stable higher string operations \cite{T2} in the context of homological conformal field theory in which homology classes of moduli spaces of Riemann surfaces give rise to string operations \cite{G}. 

As a consequence of Theorem B, we obtain the
following result on torsion elements proved in \S3. 
Let $\iota: \Omega M @>>> LM$ be
the inclusion map from the based loop space to the free loop
space. Recall that the transfer map $\iota_!: H_{*+d}(LM) @>>>
H_*(\Omega M)$ obtained by intersecting cycles with $\Omega M$ is an
algebra map with respect to the loop product in $H_*(LM)$ and the
Pontrjagin product in $H_*(\Omega M)$.

\proclaim{Theorem C} Let $M$ be an even dimensional manifold with
$\chi(M)\not=0$. Consider the following composition map 
$$
\iota_*\circ\iota_!: H_{p+d}(LM) @>{\iota_!}>> 
H_p(\Omega M) @>{\iota_*}>> H_p(LM).
$$
If $p\not=0$, then the image of $\iota_*\iota_!$ consists of
torsion elements of order a divisor of $\chi(M)$. Namely, 
$$
\chi(M)\iota_*\iota_!(a)=\chi(M)[c_0]\cdot a=0 \quad \text{\ if \ }|a|\not=d 
\text{\ for\ } a\in H_*(LM). 
$$
Thus, rationally, the composition is a trivial map if $p\not=0$.
\endproclaim

See Example 3.6 for explicit examples of this fact 
when $M$ is $S^{2n}$ or $\Bbb CP^n$.  

Since Theorem A can be quickly proved from Theorem B, we give its
proof here in the remainder of this introduction. 

\demo{Proof of Theorem A from Theorem B} Let $S(p,q)$ be a genus $0$
surface with $p$ incoming and $q$ outgoing parametrized boundary
circles, and let $T$ be a torus with $1$ incoming and $1$ outgoing
parametrized boundary circles. Then any surface $\Sigma$ of genus $g$
with $p$ incoming boundary circles and $q$ outgoing boundary circles
can be decomposed as $S(p,1)\#T\#\cdots\# T\#S(1,q)$, where $T$ appears
$g$ times. Correspondingly, the associated string operation
$\mu_{\Sigma}$ can be decomposed as
$$
\mu_{\Sigma}=\mu_{S(1,q)}\circ \mu_{T}\circ\cdots\circ \mu_{T}\circ
\mu_{S(p,1)}.
$$
Assume $g\ge1$. We compute $\mu_{T}$ using a decomposition of $T$ into
two pairs of pants corresponding to the loop coproduct and the loop
product. For any $a\in H_*(LM)$,
$$
\mu_T(a)=\mu\circ\Psi(a)
=\mu\bigl(\chi(M)[c_0]\otimes([c_0]\cdot a)\bigr)
=(-1)^d\chi(M)([c_0]\cdot[c_0])\cdot a.
$$
Since $[c_0]\cdot[c_0]=0\in H_{-d}(LM)$ by dimensional reason, we have
$\mu_T(a)=0$ for all $a\in H_*(LM)$. In view of the above decomposion
of $\mu_{\Sigma}$, this proves the vanishing of string operations
associated to surfaces of genus $g\ge1$. 

Next we assume $q\ge3$. Then $\mu_{S(1,q)}=(\mu_{S(1,q-2)}\otimes
1\otimes 1)\circ (\Psi\otimes1)\circ \Psi$. For any $a\in H_*(LM)$,
$$
(\Psi\otimes 1)\circ\Psi(a)=(\Psi\otimes
1)(\chi(M)[c_0]\otimes[c_0]\cdot a)
=\chi(M)\Psi([c_0])\otimes [c_0]\cdot a=0,
$$
since $\Psi([c_0])=0\in H_{-d}(LM\times LM)$ by dimensional
reason. Hence $\mu_{S(1,q)}=0$ for $q\ge3$.  Again, in view of the
above decomposition of $\mu_{\Sigma}$, this proves $q\ge3$ case of
Theorem A.
\qed
\enddemo

In \S2, we discuss general properties of the loop coproduct in detail and prove
Frobenius compatibility (Theorem 2.2), a symmetry property (Proposition 2.3), and coderivation property of certain
cap products (Theorem 2.4). In \S3. we prove Theorem B and related results in Theorem 3.1, and deduce their consequences including Theorem C proved in Corollary 3.3 and Corollary 3.4. We also discuss torsion properties of certin loop bracket elements in Corollary 3.5, and other miscellaneous properties of image elements of the loop coproduct in Propositions 3.7 and 3.8.  
All homology groups in this paper have integer coefficients.

\specialhead
\S 2. \quad The loop coproduct and its Frobenius compatibility
\endspecialhead

As before, let $LM$ be the free loop space of continuous maps from the
circle $S^1=\Bbb R/\Bbb Z$ to a connected oriented closed smooth
$d$-manifold $M$. Cohen and Jones \cite{CJ} gave a homotopy theoretic
description of the loop product. The loop coproduct can be described
in a similar way, and we study its properties in this section. A
description of the loop coproduct using transversal chains is
given in \cite{S}.

Let $p, p': LM @>>> M$ be evaluation maps given by
$p(\gamma)=\gamma(0)$ and $p'(\gamma)=\gamma(\frac12)$ for $\gamma\in
LM$. We consider the following diagram where $SM=(p,p')^{-1}
\bigl(\phi(M)\bigr)$ consists of loops $\gamma$ such that
$\gamma(0)=\gamma(\frac12)$, and $q$ is the restriction of $(p,p')$ to
this subspace. Let $\iota: SM @>>> LM$ be the inclusion map and let
$j:SM @>>> LM\times LM$ be given by
$j(\gamma)=(\gamma_{[0,\frac12]},\gamma_{[\frac12,1]})$. The map
$\phi: M @>>> M\times M$ is the diagonal map.
$$
\CD
LM @<{\iota}<< SM @>{j}>> LM\times LM \\
@V{(p,p')}VV @V{q}VV   @. \\
M\times M @<{\phi}<< M @. 
\endCD
$$
Then the coproduct map $\Psi$ is defined by the following composition
of maps: 
$$
\Psi=j_*\circ\iota_!: H_{*+d}(LM) @>{\iota_!}>> H_*(LM) @>{j_*}>>
H_*(LM\times LM),
$$
where $\iota_!$ is the transfer map, also called a push-forward map,
defined in the following way. Let $\pi: \nu @>>> \phi(M)$ be the
normal bundle to $\phi(M)$ in $M\times M$ and we orient $\nu$ so that
we have an oriented isomorphism $\nu\oplus T\phi(M)\cong T(M\times
M)|_{\phi(M)}$. Let $N$ be a closed tubular neighborhood of $\phi(M)$
such that $D(\nu)\cong N$, where $D(\nu)$ is the closed disc
bundle. Let $c: M\times M @>>> N/\partial N$ be the Thom collapse
map. We have the following commutative diagram: 
$$
\CD 
H^d(M\times M, M\times M-\phi(M)) @>>> H^*(M\times M) \\
@V{\cong}V{\text{excision}}V   @A{c^*}AA \\
H^d(N, N-\phi(M)) @>{\cong}>> H^d(N,\partial
N)\cong\tilde{H}^d(N/\partial N).
\endCD
$$
Let $u'\in\tilde{H}(N/\partial N)$ be the Thom class of the normal
bundle $\nu$.  Let $u''\in H^d\bigl(M\times M,M\times M-\phi(M)\bigr)$
and $u\in H^d(M\times M)$ be corresponding Thom classes. The class $u$
is characterized by the property $u\cap[M\times M]=\phi_*([M])$. Since
$u$ comes from $u''$, it is represented by a cocycle $f$ which vanish
on simplices in $M\times M$ which do not intersect with $\phi(M)$.

Let $\tilde{N}=(p,p')^{-1}(N)$ be a tubular neighborhood of $SM$ in
$LM$, and let $\tilde{c}: LM @>>> \tilde{N}/\partial\tilde{N}$ be the
Thom collapse map. Let $\tilde{u}'\in
\tilde{H}^d(\tilde{N}/\partial\tilde{N})$ and $\tilde{u}\in H^d(LM)$
be pull-backs of corresponding classes. We have
$\tilde{u}=\tilde{c}^*(\tilde{u}')$. Let $\tilde{\pi}:\tilde{N} @>>>
SM\subset LM\times LM$ be a projection map corresponding to $\pi$, and is
given as follows.  Suppose $\gamma\in\tilde{N}$ is such that
$(p,p')(\gamma)=(x_1,x_2)\in N$.  Let $\eta(t)=(\eta_1(t),\eta_2(t))$
be a path in $N$ from $(x_1,x_2)$ to $\pi(x_1,x_2)=(y,y)\in\phi(M)$
corresponding to the straight ray in the bundle $\nu$. Then
$\tilde{\pi}(\gamma)=(\eta_1^{-1}\cdot\gamma_{[0,\frac12]}\cdot\eta_2)\cdot
(\eta_2^{-1}\cdot \gamma_{[\frac12,1]}\cdot\eta_1)\in SM$. From this
description, it is obvious that $\tilde{\pi}$ is a deformation
retraction. The transfer map $\iota_!$ is defined by the following
composition of maps:
$$
\iota_!: \tilde{H}_{*+d}(LM) @>{\tilde{c}_*}>> 
\tilde{H}_{*+d}(\tilde{N}/\partial\tilde{N})
@>{\tilde{u}'\cap(\cdot)}>> H_*(\tilde{N})
@>{\tilde{\pi}_*}>{\cong}> H_*(SM).
$$
Let $s: M @>>> LM$ be the constant loop map given by $s(x)=c_x$, where
$c_x$ is the constant loop at $x\in M$. Since $p\circ s=1_M$, we have
$s^*\circ p^*=1$. The transfer map $\iota_!$ has the following
properties.

\proclaim{Proposition 2.1} \rom{(1)} The cohomology class 
$\tilde{u}\in H^d(LM)$ is given by $\tilde{u}=p^*(e_M)$, where $e_M\in
H^d(M)$ is the Euler class of $M$.

\rom{(2)} For any element $a\in H_*(LM)$, 
$$
\iota_*\iota_!(a)=p^*(e_M)\cap a.
\tag2-1
$$
In particular, $\iota_*\iota_!(s_*([M]))=\chi(M)[c_0]$, where $c_0$ is
the constant loop at the base point $x_0$ in $M$, and $\chi(M)$ is the
Euler characteristic of $M$.

\rom{(3)} For any $\alpha\in H^*(LM)$ and $b\in H_*(LM)$, 
$$
\iota_!(\alpha\cap b)=(-1)^{d|\alpha|}\iota^*(\alpha)\cap\iota_!(b).
\tag2-2
$$
\endproclaim
\demo{Proof} (1) Since the map $(p,p'):LM @>>> M\times M$ can be
factored as $LM @>{\phi}>> LM\times LM @>{p\times p'}>> M\times M$ and
$p$ and $p'$ are homotopic, we have $\tilde{u}=(p,p')^*(u)=\phi^*\circ
(p\times p')^*(u)=\phi^*\circ(p\times p)^*(u)=p^*\circ
\phi^*(u)$. Since $\phi^*(u)$ is, by definition, $(-1)^d$ times the Euler
class $e_M$ of $M$ and the Euler class is of order 2 when $d$ is odd,
we have $(-1)^de_M=e_M$. So we have $\tilde{u}=p^*(e^M)$.

(2) Although we can use a certain commutative diagram for a proof (see
    below), we first do a chain argument here in the spirit of
    \cite{CS} and \cite{S}. By barycentric subdivisions on the cycle
    $\xi$ representing $a\in H_*(LM)$, we may assume that every
    simplex of $\xi$ intersecting with $SM$ is contained in
    $\text{Int}\,(\tilde{N})$. Since cohomology classes $u\in
    H^d(M\times M)$ and $u'\in \tilde{H}^d(N/\partial N)$ come from
    the class $u''$ in $H^d\bigl(M\times M, M\times M-\phi(M)\bigr)$,
    they can be represented by cocycles $f$ and $f'$ so that $f$
    vanishes on simplices in $M\times M$ not intersecting $\phi(M)$,
    and $f'$ vanishes on simplices in $N/\partial N$ not intersecting
    $\phi(M)$. So the cocycle $\tilde{f}'=(p,p')^{\#}(f')$
    representing $\tilde{u}'$ vanishes on simplices in $N$ not
    intersecting with $SM$. Similarly, the cocycle
    $\tilde{f}=(p,p')^{\#}(f)=\tilde{c}^{\#}(\tilde{f}')$ representing
    $\tilde{u}=\tilde{c}^*(\tilde{u}')$ vanishes on simplices in $LM$
    not intersecting with $SM$, and has the same values as
    $\tilde{f}'$ on simplices in $\tilde{N}$ intersecting with
    $SM$. Since the cycle $\xi$ is fine enough, the cycles
    $\tilde{f}\cap\xi$ and $\tilde{f}'\cap \tilde{c}^{\#}(\xi)$
    representing $\tilde{u}\cap a$ and $\tilde{u}'\cap \tilde{c}^*(a)$
    are in fact identical. Since
    $\iota_!(a)=\tilde{\pi}_*\bigl(\tilde{u}'\cap
    \tilde{c}_*(a)\bigr)$ is represented by a cycle
    $\tilde{\pi}_{\#}(\tilde{f}\cap \xi)$, and $\tilde{\pi}$ is a
    deformation retraction, the two cycles
    $\tilde{\pi}_{\#}(\tilde{f}\cap \xi)$ and $\tilde{f}\cap \xi$ are
    homologous inside of $\text{Int}\,\tilde{N}$. Thus,
    $\iota_*\iota_!(a)=[\tilde{\pi}_{\#}(\tilde{f}\cap \xi)]$ and
    $\tilde{u}\cap a=[\tilde{f}\cap \xi]$ represent the same homology
    class. Hence $\iota_*\iota_!(a)=\tilde{u}\cap a=p^*(e_M)\cap a$,
    by (1).

We also give a homological proof, using the following commutative
diagram. 
$$
\CD
H_*(LM) @>{\tilde{c}_*}>> H_*(\tilde{N}, \partial\tilde{N}) 
@>{\tilde{u}'\cap(\ \cdot\ )}>> 
H_{*-d}(\tilde{N}) @>{\tilde{\pi}_*}>{\cong}>  H_{*-d}(SM) \\
@|   @V{(\iota_N )_*}V{\cong}V  @V{(\iota_N )_*}VV  @V{\iota_*}VV  \\
H_*(LM) @>{j_*}>>  H_*(LM, LM-SM)  @>{\tilde{u}''\cap(\ \cdot\ )}>> 
H_{*-d}(LM) @= H_{*-d}(LM)
\endCD
$$
where $\iota_N: \tilde{N} @>>> LM$ is an inclusion map. 
Here, the class $\tilde{u}''$ is given by $\tilde{u}''=(p,p')^*(u'')$,
and it satisfies $\tilde{u}'=\iota_N^*(\tilde{u}'')$. 
Thus, for $a\in H_*(LM)$, the
commutative diagram shows $\iota_*\iota_!(a)=\tilde{u}''\cap j_*(a)
=j^*(\tilde{u}'')\cap a=\tilde{u}\cap a$. The above chain argument
gives geometric meaning to the commutative diagram above. 

When $a=s_*([M])$, we have $\iota_*\iota_!\bigl(s_*([M])\bigr)
    =p^*(e_M)\cap s_*([M])=s_*\bigl(s^*p^*(e_M)\cap[M]\bigr)$. Since
    $p\circ s=1$, this is equal to
    $s_*\bigl(\chi(M)[x_0]\bigr)=\chi(M)[c_0]$.

(3) We compute. By definition of $\iota_!$, we have
$$
\iota_!(\alpha\cap b)=\tilde{\pi}_*\bigl(\tilde{u}\cap \tilde{c}_*
    (\alpha\cap b)\bigr)=
\tilde{\pi}_*\bigl(\tilde{u}'\cap
(\iota_N^*(\alpha)\cap\tilde{c}_*(b)\bigr)\bigr) 
=(-1)^{|\alpha|d}\tilde{\pi}_*\bigl(\iota_N^*(\alpha)\cap 
\bigl(\tilde{u}'\cap\tilde{c}_*(b)\bigr)\bigr).
$$
Since $\iota_N^*(\alpha)=\tilde{\pi}^*\iota^*(\alpha)$, the last 
formula becomes $\iota^*(\alpha)\cap \tilde{\pi}_*
\bigl(\tilde{u}'\cap\tilde{c}_*(b)\bigr)
=\iota^*(\alpha)\cap\iota_!(b)$, times the sign. 
This completes the proof. 
\qed
\enddemo

Next we recall a homotopy theoretic description of the loop product
from \cite{CJ}. We consider the following diagram, where
$LM\times_MLM$ denotes the set $(p\times
p)^{-1}\bigl(\phi(M)\bigr)$ consisting of pairs $(\gamma, \eta)$ of loops 
such that $\gamma(0)=\eta(0)$, and $\iota(\gamma,\eta)$ denotes the usual loop
multiplication $\gamma\cdot\eta$.
$$
\CD
LM\times LM @<{j}<< LM\times_MLM @>{\iota}>> LM \\
@V{p\times p}VV @V{q}VV    @.\\
M\times M @<{\phi}<< M   @.
\endCD
$$
Then for $a,b\in H_*(LM)$, the loop product $a\cdot b$ is defined by 
$$
a\cdot b=(-1)^{d(|a|-d)}\iota_*j_!(a\times b).
$$
Here, as before, the transfer map $j_!$ is defined using the Thom
class $u'\in \tilde{H}^d(N/\partial N)$ and its pull back to the
tubular neighborhood $\tilde{N}=(p\times p)^{-1}(N)$. The sign
$(-1)^{d(|\alpha|-d)}$ is natural since on the right hand side, the
map $j_!$, which represents the content of the loop product, is in
front of $a$, whereas on the left hand side, the dot representing the loop
product is between $a$ and $b$. Switching the order of $j_!$ and $a$
yields the sign $(-1)^{d|\alpha|}$. The sign $(-1)^d$ comes from our
choice of the orientation of the normal bundle $\nu$ so that $[M]\in
H_d(LM)$ acts as the unit. Note that the $|a|-d$ is the degree of $a$
in the loop algebra $\Bbb H_*(LM)=H_{*+d}(LM)$.

For further discussion, we need transfer maps defined in the following
general context. Let $\iota: K @>>> M$ be a smooth embedding of
oriented closed smooth manifolds and let $\nu$ be its normal bundle
oriented by $\nu\oplus\iota_*(TK)\cong TM|_{\iota(K)}$. Let $u'$ be
the Thom class of $\nu$ and let $u\in H^{d-k}(M)$ be the corresponding
Thom class for the embedding $\iota$, where $d$ and $k$ are dimensions
of $M$ and $K$. With the above choice of the orientation on $\nu$, we
have $u\cap [M]=\iota_*([K])$, which characterizes the class $u$. Had
we oriented $\nu$ by $\iota_*(TK)\oplus\nu\cong TM|_{\iota(K)}$, then
we would have obtained $u\cap[M]=(-1)^{k(d-k)}\iota_*([K])$.

Let $p: E @>>> M$ be a Hurewicz fibration, and let $E_K$ be its
pull-back over $K$ via the embedding $\iota$. Let $\iota: E_K @>>> E$
be the inclusion of fibrations. Proceeding as before, we can define a
transfer map.
$$
\iota_! : H_{*+d}(E) @>>> H_{*+k}(E_K), \text{\ \ such that \ \ }
\iota_*\iota_!(a)=p^*(u)\cap a \text{\ \ for any \ }a\in H_*(E).
$$
We remark that with the above choice of the orientation on the 
normal bundle $\nu$, the transfer map between base manifolds satisfies 
$\iota_!([M])=[K]$. Also, it can be verified that for a composition of 
smooth embeddings $K @>{g}>> L @>{f}>> M$ and the associated induced 
inclusions of fibrations $E_K @>{g}>> E_L @>{f}>> E$, we have 
$(f\circ g)_!=g_!\circ f_!$. 

The loop product enjoys the Frobenius compatibility with respect to the loop
coproduct, in the following sense. This is discussed in \cite{S} from
the point of view of chains. Here, we give a homotopy theoretic proof 
with precise determination of signs.  

For $a\in H_*(LM)$ and $c\in H_*(LM\times LM)$, let $a\cdot c$ be
defined by $(\iota\times 1)_*\circ(j\times 1)_!(a\times
c)=(-1)^{d(|a|-d)}a\cdot c$ using the following diagram
$$
\CD
(LM\times LM)\times LM  @>{j\times 1}>> 
(LM\times_M LM)\times LM @>{\iota\times 1}>> 
LM\times LM \\
@V{p_1\times p_2}VV  @V{p_1}VV  \\
M\times M @<{\phi}<< M @. 
\endCD
$$
where $p_1\times p_2$ denotes projections from the first and second
factor. If $c$ is of the form of a cross product $b\times c$, then
$a\cdot(b\times c)=(a\cdot b)\times c$. Similarly, an element $c\cdot
a$ is defined by $(1\times\iota)_*(1\times j)_!(c\times
a)=(-1)^{d(|c|-d)}c\cdot a$ using a similar diagram.

\proclaim{Theorem 2.2} The loop product and the loop coproduct satisfy 
Frobenius compatibility, namely, for $a,b\in H_*(LM)$, 
$$
\Psi(a\cdot b)=(-1)^{d(|a|-d)}a\cdot\Psi(b)=\Psi(a)\cdot b.
\tag2-3
$$
\endproclaim
\demo{Proof} For convenience, we introduce a space $L^rM$ of
continuous loops from a circle of length $r>0$ to $M$. We let
$L'M=L^{\frac13}M$ and $L''M=L^{\frac23}M$. We identify $SM\subset
L^{2r}M$ with $L^rM\times_M L^rM$. We have the following
commutative diagram of inclusions: 
$$
\CD
L'M\times L''M @<{1\times \iota}<< 
L'M\times L'M\underset{M}\to{\times} L'M 
@>{1\times j}>> L'M\times L'M\times L'M \\
@A{j}AA @A{j_1=(j\times_M1)}AA   @A{j\times 1}AA \\
L'M\underset{M}\to{\times}L''M @<{\iota_1=(1\times_M\iota)}<< 
L'M\underset{M}\to{\times}L'M\underset{M}\to{\times}L'M 
@>{j_2=(1\times_Mj)}>> L'M\underset{M}\to{\times}L'M\times L'M  \\
@V{\iota}VV  @V{\iota_2=(\iota\times_M1)}VV    @V{\iota\times 1}VV   \\
LM @<{\iota}<< L''M\underset{M}\to{\times}L'M @>{j}>> L''M\times L'M.
\endCD
$$
The base manifolds of fibrations in the above diagram form the
following diagram which we use to compute Thom classes of embeddings,
which in turn are used to construct transfer maps.
$$
\CD
M\times M\times M @<{1\times \phi}<< M\times M @>{1\times \phi}>> 
M\times M\times M \\
@A{\phi\times 1}AA     @A{\phi}AA    @A{\phi\times 1}AA   \\
M\times M   @<{\phi}<<  M @>{\phi}>>  M\times M  \\
@V{\phi\times 1}VV  @V{\phi}VV @V{\phi\times 1}VV   \\
M\times M\times M  @<{\phi_{13}}<<    M\times M 
@>{\phi_{13}}>>  M\times M\times M
\endCD
$$
where $\phi_{13}(x,y)=(x,y,x)$, or $\phi_{13}=(1\times T)(\phi\times 1)$
and $T:M\times M @>>> M\times M$ is the switching map. Here, for example, 
the fibration $p: L''M\times_M L'M @>>> M\times M$ is given by 
$p(\gamma,\eta)=(\gamma(0)=\eta(0), \gamma(\frac13))$, and the fibration
$p : L'M\times_M L''M @>>> M\times M$ is given by $p(\gamma, \eta)=
(\gamma(0)=\eta(0),\eta(\frac13))$. 

To prove $\Psi(a\cdot b)=(-1)^{d(|a|-d)}a\cdot\Psi(b)$, we examine the
following induced homology diagram with transfers in which we replaced
$L'M$ and $L''M$ by their homeomorphic copy $LM$. 
$$
\CD
H_*(LM\times LM) @>{(1\times \iota)_!}>{=(-1)^d1\times \iota_!}> 
H_{*-d}(LM\times LM\underset{M}\to{\times} LM)  @>{(1\times j)_*}>> 
H_{*-d}(LM\times LM\times LM) \\
@V{\tilde{j}_!=j_!}VV   @V{(j_1)_!}VV   @V{(j\times 1)_!=j_!\times 1}VV
\\
H_{*-d}(LM\underset{M}\to{\times}LM) @>{(\iota_1)_!}>> 
H_{*-2d}(LM\underset{M}\to{\times}LM\underset{M}\to{\times}LM) @>{(j_2)_*}>> 
H_{*-2d}(LM\underset{M}\to{\times}LM\times LM)  \\
@V{\iota_*}VV   @V{(\iota_2)_*}VV   @V{(\iota\times 1)_*}VV   \\
H_{*-d}(LM)  @>{\tilde{\iota}_!=(-1)^d\iota_!}>> 
H_{*-2d}(LM\underset{M}\to{\times}LM) 
@>{j_*}>> H_{*-2d}(LM\times LM)
\endCD
$$
In the above, the transfer maps $\tilde{j}_!, \tilde{\iota}_!$
indicate that Thom classes used to define these transfer maps may be
different in signs from Thom classes used to define transfers
$\iota_!$ and $j_!$.

The top left square and the bottom right square commute because of the
functorial properties of transfer maps and induced maps. We
examine the commutativity of the bottom left square. Since the
corresponding square of fibrations commutes, the homology square
with induced maps and transfer maps
commutes up to a sign. To determine this sign, for $a\in
H_*(LM\times_M LM)$, we compare $\iota_*(\iota_2)_*(\iota_1)_!(a)$ and
$\iota_*\tilde{\iota}_!\iota_*(a)$ in $H_*(LM)$. Let $u\in H^d(M\times
M)$ be the Thom class for the embedding $\phi: M @>>> M\times M$. Then
the Thom class for the embedding $\phi_{13} : M\times M @>>>
M\times M\times M$ is given by $(-1)^du_{13}$, where 
$u_{13}=(1\times T)^*(u\times 1)=\sum_i (u_i'\times 1\times u_i'')$ if 
$u=\sum_i u_i'\times u_i''$. Hence
$(\iota_*\tilde{\iota}_!)\iota_*(a)=(-1)^dp^*(u_{13})\cap
\iota_*(a)$, where the map $p: LM @>>> M\times M\times M$ is a
fibration given by $p(\gamma)=(\gamma(0),
\gamma(\frac13),\gamma(\frac23))$. On the other hand, using the
commutativity of the induced homology square, we have
$\iota_*(\iota_2)_*(\iota_1)_!(a)=\iota_*(\iota_1)_*(\iota_1)_!(a)
=\iota_*(p^*(u)\cap a)$, since the Thom class for the embedding
$\iota_1$ is $p^*(u)$. Since $u=(\phi\times 1)^*(u_{13})$, we have 
$p^*(u)=p^*((\phi\times 1)^*(u_{13}))=\iota^*(p^*(u_{13}))$. 
Hence $\iota_*(p^*(u)\cap a)=p^*(u_{13})\cap
\iota_*(a)$. Collecting our computations, we have that
$\iota_*(\iota_2)_*(\iota_1)_!(a)=p^*(u_{13})\cap
\iota_*(a)$. Comparing with the formula above for
$\iota_*\tilde{\iota}_!\iota_*(a)$, we see that the sign difference
between $(\iota_2)_*(\iota_1)_!(a)$ and $\tilde{\iota}_!\iota_*(a)$
is given by $(-1)^d$. Hence the square commutes up to $(-1)^d$. 

Similar argument shows that the top right square in the homology
diagram actually commutes. 

Next we examine transfer maps in the diagram. For the top horizontal left
transfer $(1\times\iota)_!$, since the Thom class of the embedding
$1\times\phi : M\times M @>>> M\times M\times M$ is $(-1)^d(1\times
u)$, 
$$
\multline
(1\times\iota)_*(1\times \iota)_!(a\times b)
=(-1)^dp^*(1\times u)\cap(a\times b) 
=(-1)^{d+d|a|}a\times(p^*(u)\cap b)  \\
=(-1)^{d+d|a|}a\times\iota_*\iota_!(b)
=(-1)^d(1\times \iota)_*(1\times(\iota)_!)(a\times b).
\endmultline
$$  
for $a,b\in H_*(LM)$, Thus, $(1\times \iota)_!=(-1)^d1\times
(\iota)_!$, as indicated in the diagram. Similarly, we can verify that
for the vertical top right transfer map, we have $(j\times
1)_!=j_!\times 1$. For the vertical top left transfer $\tilde{j}_!$
associated to the Thom class for the embedding $\phi\times 1: M\times
M @>>> M\times M\times M$ coincides with the transfer $j_!$ associated
to the Thom class for the embedding $\phi: M @>>> M\times M$. The
bottom left horizontal transfer map $\tilde{\iota}_!$ associated to
the Thom class $(-1)^du_{13}$ for the embedding $\phi_{13} : M\times M
@>>> M\times M\times M$ coincides with $(-1)^d\iota_!$, where
$\iota_!$ is the transfer associated to the Thom class $u$ of the
embedding $\phi: M @>>> M\times M$.

Hence for $a,b\in H_*(LM)$, tracing the diagram from the top left corner to 
the bottom right corner via bottom left corner, we get 
$$
j_*(\tilde{\iota})_!\iota_*(\tilde{j})_!(a\times b)
=j_*(-1)^d\iota_!\bigl((-1)^{d(|a|-d)}a\cdot b\bigr)
=(-1)^{d+d(|a|-d)}\Psi(a\cdot b).
$$
Following the diagram via the top right corner, we get 
$$
\multline
(\iota\times 1)_*(j\times 1)_!(1\times j)_*(1\times\iota)_!(a\times b)
=(\iota_*\times 1)(j_!\times 1)(1\times j_*)(-1)^d(1\times\iota_!)(a\times b) \\
=(-1)^{d+|a|d}(\iota_*j_!\times 1)(a\times\Psi(b))
=(-1)^{d+|a|d+d(|a|-d)}a\cdot\Psi(b).
\endmultline
$$
Since the entire diagram commutes up to $(-1)^d$, we finally get 
$\Psi(a\cdot b)=(-1)^{d(|a|-d)}a\cdot\Psi(b)$. 

To prove the other identity $\Psi(a\cdot b)=\Psi(a)\cdot b$, 
we consider the induced homology diagram with transfers flowing from the 
bottom right corner to the top left corner given below. 
$$
\CD
H_{*-2d}(LM\times LM) @<{(1\times \iota)_*}<< 
H_{*-2d}(LM\times LM\underset{M}\to{\times} LM)  
@<{(1\times j)_!}<{=(-1)^d(1\times j_!)}< 
H_{*-d}(LM\times LM\times LM) \\
@A{\tilde{j}*}AA   @A{(j_1)*}AA   @A{(j\times 1)_*}AA
\\
H_{*-2d}(LM\underset{M}\to{\times} LM) @<{(\iota_1)_*}<<
H_{*-2d}(LM\underset{M}\to{\times} LM\underset{M}\to{\times}LM) 
@<{(j_2)!}<< 
H_{*-d}(LM\underset{M}\to{\times} LM\times LM)  \\
@A{\tilde{\iota}_!=\iota_!}AA   @A{(\iota_2)_!}AA  
 @A{(\iota\times 1)_!=\iota_!\times 1}AA  \\
H_{*-d}(LM)  @<{{\iota}_*}<< H_{*-d}(LM\underset{M}\to{\times}LM) 
@<{\tilde{j}_!=(-1)^dj_!}<< H_*(LM\times LM)
\endCD
$$
where the transfer maps along the perimeter has been identified as shown. 
Using similar methods, all the squares commute except the top right 
one which commutes up to $(-1)^d$. With this information, 
following the diagram via top right corner gives $(-1)^{d+d|a|}\Psi(a)\cdot b$, 
and following the diagram via the bottom left corner gives 
$(-1)^{d+d(|a|-d)}\Psi(a\cdot b)$. Since the entire diagram commutes 
up to $(-1)^d$, we obtain the identity $\Psi(a\cdot b)=\Psi(a)\cdot b$. 
This completes the proof. 
\qed
\enddemo

Note that in the same diagram of fibrations, if we consider an induced
homology diagram with transfers flowing from the top right corner to
the bottom left corner, or a diagram flowing from the bottom left
corner to the top right corner, we obtain homotopy theoretic proofs of
associativity of the loop product \cite{CJ} and the coassociativity of
the loop coproduct.

Next we show that $\Psi$ is symmetric. Let $T:LM\times LM @>>>
LM\times LM$ be the switching map. 

\proclaim{Proposition 2.3} The loop coproduct is symmetric in the
sense that 
$$
T_*\bigl(\Psi(a)\bigr)=\Psi(a)
$$ 
for any $a\in H_*(LM)$. 
\endproclaim
\demo{Proof} We consider the following commutative diagram: 
$$
\CD 
LM @<{\iota}<< LM\times_M LM @>{j}>> LM\times LM \\
@V{R_{\frac12}}VV   @V{T}VV   @V{T}VV    \\
LM @<{\iota}<< LM\times _M LM @>{j}>> LM\times LM
\endCD
$$
Here, as before, we identify $SM$ with $LM\times_M LM$, and
$R_{\frac12}$ is the rotation of loops by $\frac12$, that is
$R_{\frac12}(\gamma)(t)=\gamma(t+\frac12)$. The left square commutes
because
$R_{\frac12}\circ\iota(\gamma,\eta)=R_{\frac12}(\gamma\cdot\eta)
=\eta\cdot\gamma=\iota\circ T(\gamma,\eta)$.  The Thom class for the
embedding $\iota$ is given by $\tilde{u}=p^*(e_M)$. Since
$R_{\frac12}\simeq 1$, we have
$R^*_{\frac12}(\tilde{u})=\tilde{u}$. Thus the Thom classes for two
$\iota$'s are compatible and we have $T_*\circ \iota_!=\iota_!\circ
{R_{\frac12}}_*=\iota_!$. Thus the above commutative diagram implies
$T_*\bigl(\Psi(a)\bigr)=T_*\circ j_*\circ \iota_!(a)=j_*\circ
T_*\circ\iota_!(a)=j_*\circ\iota_!(a)=\Psi(a)$.
\qed
\enddemo

The loop coproduct behaves well with respect to cap products with
cohomology classes in $H^*(LM)$ arising from $\alpha\in H^*(M)$. Let
$p: LM @>>> M$ be the base point map. For the evaluation map $e=p\circ
\Delta: S^1\times LM @>>> M$, let $e^*(\alpha)=1\times
p^*(\alpha) +\{S^1\}\times\Delta\bigl(p^*(\alpha)\bigr)$, where $\{S^1\}$ 
is the fundamental cohomology class for $S^1$.

\proclaim{Theorem 2.4} Let $\alpha\in H^*(M)$ and $b\in H_*(LM)$. 

\rom{(1)} The cap product with $p^*(\alpha)$ satisfies Frobenius
    compatibility with respect to the loop coproduct\rom{:} 
$$
\Psi\bigl(p^*(\alpha)\cap b\bigr)
=(-1)^{d|\alpha|}\bigl(p^*(\alpha)\times 1\bigr)\cap\Psi(b)
=(-1)^{d|\alpha|}\bigl(1\times p^*(\alpha)\bigr)\cap\Psi(b).
\tag2-4
$$

\rom{(2)}  The cap product with $\Delta\bigl(p^*(\alpha)\bigr)$
behaves as a coderivation with respect to the loop coproduct\rom{:}
$$
\Psi\bigl(\Delta\bigl(p^*(\alpha)\bigr)\cap b\bigr)
=(-1)^{d(|\alpha|-1)}\bigl[\Delta\bigl(p^*(\alpha)\bigr)\times 1 +
1\times \Delta\bigl(p^*(\alpha)\bigr)\bigr]\cap\Psi(b).
\tag2-5
$$
\endproclaim
\demo{Proof}  From the definition of the loop coproduct and a property
(2-2) of the transfer $\iota_!$, we have $\Psi\bigl(p^*(\alpha)\cap
b\bigr)=j_*\circ\iota_!\bigl(p^*(\alpha)\cap b\bigr)
=(-1)^{d|\alpha|}j_*\bigl(\iota^*p^*(\alpha)\cap \iota_!(b)\bigr)$.
To understand $\iota^*p^*(\alpha)$, we consider the following
commutative diagram.
$$
\CD 
LM @<{\iota}<< LM\times_M LM @>{j}>> LM\times LM  \\
@V{p\times p'}VV  @V{q}VV   @V{p\times p}VV   \\
M\times M @<{\phi}<<  M @>{\phi}>> M\times M \\
@V{\pi_1}VV  @|  @V{\pi_i}VV \\
M @= M @= M
\endCD
$$
where $p'(\gamma)=\gamma(\frac12)$, and $\pi_i$ for $i=1,2$ is the
projection onto the $i$th factor. From the diagram, we
have $\iota^*p^*(\alpha)=q^*(\alpha)=j^*(p\times p)^*\pi_i^*(\alpha)$,
which is equal, for $i=1,2$, to $j^*\bigl(p^*(\alpha)\times 1\bigr)$
and to $j^*\bigl(1\times p^*(\alpha)\bigr)$. For $i=1$ case, 
$$
(-1)^{d|\alpha|} \Psi\bigl(p^*(\alpha)\cap b\bigr)
=j_*\bigl(j^*\bigl(p^*(\alpha)\times 1)\cap \iota_!(b)\bigr)
=\bigl(p^*(\alpha)\times 1\bigr)\cap j_*\iota_!(b)
=\bigl(p^*(\alpha)\times 1\bigr)\cap \Psi(b). 
$$
Similarly, for the case $i=2$, we obtain
$(-1)^{d|\alpha|}\Psi\bigl(p^*(\alpha)\cap b\bigr) =\bigl(1\times
p^*(\alpha)\bigr)\cap \Psi(b)$.

For (2), first we note that 
$$
\Psi\bigl(\Delta(p^*(\alpha))\cap
 b\bigr)=j_*\iota_!\bigl(\Delta(p^*(\alpha))\cap b\bigr)
=(-1)^{d(|\alpha|-1)}j_*\bigl(\iota^*\Delta(p^*(\alpha))
\cap\iota_!(b)\bigr). 
$$
We need to understand $\iota^*\bigl(\Delta(p^*(\alpha))\bigr)$.  For
 this purpose, we introduce some notations. Let $I_1=[0,\frac12]$ and
 $I_2=[\frac12, 1]$. Let $r: S^1=I/\partial I @>>> I/\{0,\frac12,
 1\}=S^1_1\vee S^1_2$, where $S^1_i=I_i/\partial I_i$ for $i=1,2$, be
 an identification map. Let $\iota_i : S^1_i @>>> S^1_1\vee S^1_2$ be
 the inclusion map for $i=1,2$. We consider the following commutative
 diagram.
$$
\CD
        @. S^1_i\times(LM\times_M LM) @>{1\times j}>>
        S^1_i\times(LM\times LM)  \\
@.     @V{\iota_i\times 1}VV  @V{1\times \pi_i}VV   \\
S^1\times (LM\times_M LM) @>{r\times 1}>> 
(S^1_1\vee S^1_2)\times(LM\times_M LM)   @. 
S^1_i\times LM\cong S^1\times LM \\ 
@V{1\times\iota}VV  @V{e'}VV   @V{e}VV \\
S^1\times LM @>{e}>> M  @= M
\endCD
$$
where $e'(t,\gamma,\eta)$ is given by $\gamma(2t)$ for $0\le
t\le\frac12$, and $\eta(2t-1)$ for $\frac12\le t\le 1$. Let
${e'}^*(\alpha)=1\times\iota^*p^*(\alpha)+ \{S^1_1\}\times
\Delta_1(\alpha) +\{S^1_2\}\times\Delta_2(\alpha)$, where the first
term is due to a fact that $e'$ restricted to $\{0\}\times(LM\times_M
LM)$ is given by $p\circ\iota$. Since $e^*(\alpha)=1\times p^*(\alpha)
+\{S^1\}\times\Delta p^*(\alpha)$ and $r^*(\{S^1_i\})=\{S^1\}$ 
for $i=1,2$, the commutativity of the left
bottom square implies that 
$$
\iota^*\Delta\bigl(p^*(\alpha)\bigr) 
=\Delta_1(\alpha)+\Delta_2(\alpha).
$$
We need to identify $\Delta_i(\alpha)$ for $i=1,2$. The commutativity
of the right square implies that, for $i=1$, $(1\times j)^*(1\times
\pi_1)^*e^*(\alpha)=1\times j^*(p^*(\alpha)\times 1)+\{S^1\}\times
j^*\bigl(\Delta(p^*(\alpha))\times 1\bigr)$ is equal to 
$(\iota_1\times 1)^*{e'}^*(\alpha)=1\times\iota^*p^*(\alpha)+
\{S^1\}\times\Delta_1(\alpha)$. Hence
$\Delta_1(\alpha)=j^*\bigl(\Delta(p^*(\alpha))\times
1\bigr)$. Similarly, the $i=2$ case implies that $\Delta_2(\alpha)
=j^*\bigl(1\times\Delta(p^*(\alpha))\bigr)$. Combining the above
calculations, we have 
$$
\align
\Psi\bigl(\Delta(p^*(\alpha))\cap b\bigr)
&=(-1)^{d(|\alpha|-1)}\bigl(\iota^*\Delta(p^*(\alpha))\cap
\iota_!(b)\bigr) \\
&=(-1)^{d(|\alpha|-1)}j_*\bigl(j^*\bigl(\Delta(p^*(\alpha))\times
1+1\times \Delta(p^*(\alpha))\bigr)\cap \iota_!(b)\bigr) \\
&=(-1)^{d(|\alpha|-1)}\bigl[\Delta(p^*(\alpha))\times
1+1\times \Delta(p^*(\alpha))\bigr]\cap \Psi(b). 
\endalign
$$
This proves the coderivation property.
\qed
\enddemo

\specialhead
\S3 Properties of the loop coproduct and their consequences
\endspecialhead

So far we have proved various algebraic properties of the loop
coproduct. These properties turn out to be strong enough to force the
loop coproduct to be given by a very simple formula, given in the next
theorem. Let $s:M @>>> LM$ be the constant loop map given by
$s(x)=c_x$, where $c_x$ is the constant loop at $x\in M$. Recall that
we assume that $M$ is connected with base point $x_0$, and let
$c_0$ be the constant loop at the base point. 

The connected components of $LM$ are in 1:1 corespondence to the set of free homotopy classes of loops $[S^1.M]$, which is in 1:1 correspondence with conjugacy classes of $\pi_1(M)$. Let 
$$
LM=(LM)_{[1]}\cup\bigcup_{[\alpha]\not=[1]}(LM)_{[\alpha]},
$$
be the decomposition of $LM$ into its components, where $[\alpha]$'s are conjugacy classes in $\pi_1(M)$. 

\proclaim{Theorem 3.1} Let $M$ be a connected oriented closed
smooth $d$-manifold. 

\rom{(1)} Let $p\ge0$ and let $a_1,a_2,\dots,a_p\in H_*(LM)$. The
loop coproduct on the 
    loop product of these elements is given by the following formula, 
    for each $0\le\ell\le p$. 
$$
\Psi(a_1a_2\cdots a_p)
=\chi(M)\bigl([c_0]\cdot a_1\cdot a_2 \cdots a_{\ell}\bigr)
\otimes\bigl([c_0]\cdot a_{\ell+1} \cdots  a_p\bigr)
\in H_*(LM)\otimes H_*(LM).
\tag3-1
$$
In particular, for the unit $1=s_*([M])\in H_d(LM)=\Bbb H_0(LM)$ of the loop
homology algebra, its coproduct is given by 
$$
\Psi(1)=\chi(M)[c_0]\otimes[c_0]\in H_0(LM)\otimes H_0(LM)\cong
H_0(LM\times LM). 
\tag3-2
$$
When $p=1$, the formula for $a\in H_*(LM)$ for $\ell=0,1$ becomes 
$$
\Psi(a)=\chi(M)\bigl([c_0]\cdot a\bigr)\otimes[c_0]
=\chi(M)[c_0]\otimes \bigl([c_0]\cdot a\bigr).
\tag3-3
$$

\rom{(2)} If $|a|\not=d$, then $\Psi(a)=0$. If $|a|=d$, then $\Psi(a)=n[c_0]\otimes[c_0]$ for some $n\in\Bbb Z$.  Thus, $\text{\rm Im}\Psi=\Bbb Z[c_0]\otimes[c_0]$. 

\medskip

\rom{(3)} Suppose $a\in H_d\bigl((LM)_{[\alpha]}\bigr)$ be a degree $d$ homology class in $[\alpha]$-component of $LM$. If $[\alpha]\not=[1]$, then $\Psi(a)=0$. 

\medskip

\rom{(4)} Suppose $a\in H_d\bigl((LM)_{[1]}\bigr)$, and suppose it is of the
form $a=ks_*([M])+(\text{decomposables})$ in the loop algebra $H_*(LM)$ 
for some $k\in\Bbb Z$, then
$$
\Psi(a)=k\chi(M)[c_0]\otimes[c_0].
$$
\endproclaim
\demo{Proof} First, we prove the formula for $\Psi(1)$. Since
    $1=s_*([M])$ has degree $d$, and $\iota_!$ decreases degree by
    $d$, we have $\iota_!(1)\in H_0(LM\times_M LM)$. Since $M$ is
    connected, connected components of $LM$ are in 1:1
    correspondence with conjugacy classes of $\pi_1(M)$. Let
    $L_0M$ be the component consisting of contractible loops so that
    $c_0\in L_0M$. Note that $L_0M\times_ML_0M$ is also connected, and
    $H_0(L_0M\times_ML_0M)\cong\Bbb Z$ is generated by
    $[(c_0,c_0)]$. So we may write $\iota_!(1)=m[(c_0,c_0)]$ for some
    $m\in \Bbb Z$. Since $\iota_*:H_0(L_0M\times_ML_0M) @>>> H_0(L_0M)$
    is an isomorphism with $\iota_*([(c_0,c_0)])=[c_0]$, 
and since (2-1) implies 
    $\iota_*\iota_!(1)=p^*(e_M)\cap s_*([M])
    =s_*(e_M\cap[M])=\chi(M)[c_0]$, we have
    $\iota_!(1)=\chi(M)[(c_0,c_0)]$. Hence
    $\Psi(1)=j_*\iota_!(1)=\chi(M)[c_0]\otimes[c_0]$. 

For $a_1,a_2,\dots,a_p\in H_*(LM)$ and for $0\le\ell\le p$, 
the Frobenius compatibility (2-3) implies
$$
\align
\Psi(a_1\cdot a_2\cdots a_p)
&=(-1)^{d(|a_1|+\cdots+|a_{\ell}|-d\ell)}
(a_1\cdots a_{\ell})\cdot \Psi(1)\cdot a_{\ell+1}\cdots a_p \\
&=(-1)^{d(|a_1|+\cdots+|a_{\ell}|-d\ell)}\chi(M) a_1\cdots
a_{\ell}\cdot[c_0]\otimes[c_0]\cdot a_{\ell+1}\cdots a_p \\
&=\chi(M)\bigl([c_0]\cdot a_1\cdots a_{\ell}\bigr)\otimes
\bigl([c_0]\cdot a_{\ell+1}\cdots a_p\bigr).
\endalign
$$
Here, we used the graded commutativity in the loop homology 
algebra given by 
$$
a\cdot b=(-1)^{(|a|-d)(|b|-d)}b\cdot a, \qquad a,b\in H_*(LM).
$$
When $p=1$, we get the formula for $\Psi(a)$ given in (3-3). Note that the formula is
compatible with the symmetry formula $T_*\Psi(a)=\Psi(a)$ in Proposition 2.3. Note also
that our formula tells us that the image of $\Psi$ is contained in the tensor product 
$H_*(LM)\otimes H_*(LM)\subset H_*(LM\times LM)$, essentially because
$\Psi(1)$ is by (3-2). 

(2) From the formula (3-3),  the value $\Psi(a)$ must be an integral multiple of $[c_0]\otimes [c_0]\in H_0(LM)\otimes H_0(LM)$. Since $\Psi$ lowers degree by $d$, if $|a|\not=d$, we must have $\Psi(a)=0$. 

(3) Let $a\in H_d\bigl((LM)_{[\alpha]}\bigr)$. We show that if $\Psi(a)\not=0$, then $[\alpha]=[1]$. By (2), $\Psi(a)$ must be of the form $n[c_0]\otimes [c_0]$ for some $n\in\Bbb Z$.  Compaing with (3-3), if $\Psi(a)\not=0$, then $[c_0]\cdot a=k[c_0]$ for some $k\not=0$, which is a homology class of finite union of contractible loops.
Thus $a$ must be represented by a cycle in the space of contractible loops $(LM)_{[1]}$. Hence we have $[\alpha]=[1]$.  

(4) By (2), if $|a|\not=d$, we must have $\Psi(a)=0$, which is equivalent to 
$$
\chi(M)[c_0]\cdot a=0, \qquad a\in H_*(LM) \text{\ with\ }|a|\not=d.
\tag3-4
$$
Now suppose $|a|=d$ and $a$ is decomposable of the form
$a=\sum_ib_i'\cdot b_i''$ with $|b_i'|\not=d$ for all $i$, then 
$\Psi(a)=\sum_i\chi(M)[c_0]\cdot b_i'\otimes[c_0]\cdot b_i''=0$ by (3-1).
 Thus, if $a$ is of the form
$a=ks_*([M])+(\text{decomposables})$, then
$\Psi(a)=\Psi\bigl(ks_*([M])\bigr)=k\chi(M)[c_0]\otimes[c_0]$. 
\qed
\enddemo

Implications of Theorem 3.1 are rather striking. First, we
start with straightforward corollaries whose proofs are obvious. 

\proclaim{Corollary 3.2} Let $M$ be a connected closed oriented smooth
manifold. If its Euler characteristic is zero, then the loop coproduct
vanishes identically. 

In particular, if $M$ is odd dimensional, then the loop coproduct
vanishes identically. 
\endproclaim

For example, the loop coproduct vanishes in $H_*(LS^{2n+1})$. The
above Corollary 3.2 was also observed in \cite{S}. 

Next, we examine torsion elements in loop homology.

\proclaim{Corollary 3.3} Assume that $\chi(M)\not=0$ for a connected
closed oriented smooth $d$-manifold $M$. For any element $a\in H_*(LM)$ with
$|a|\not=d$, the element $[c_0]\cdot a$ is either $0$ or a torsion
element of order a divisor of $\chi(M)$. 
\endproclaim
\demo{Proof} In the proof of Theorem 3.1, we noted that
$\chi(M)[c_0]\cdot a=0$ if $|a|\not=d$ in (3-4). Since $\chi(M)\not=0$, the
conclusion follows. 
\qed
\enddemo

When $|a|=d$, the element $[c_0]\cdot a$ lies in $H_0(LM)$, so it is
either $0$ or torsion free. 

Let $\iota: \Omega M @>>> LM$ be the inclusion map from the based loop
space to the free loop space. Recall that we have an algebra map 
$$
\iota_!: H_{*+d}(LM) @>>> H_*(\Omega M)
$$
from the loop algebra to the Pontrjagin ring, where $d=\dim M$.  

\proclaim{Corollary 3.4} Suppose $\chi(M)\not=0$ for a closed oriented
smooth $d$-manifold $M$. Then for $p\not=0$, the image of the composition 
$$
\iota_*\circ \iota_!: H_{p+d}(LM) @>>> H_p(\Omega M) @>>> H_p(LM)
$$
consists entirely of torsion elements of order a divisor of $\chi(M)$. 
\endproclaim
\demo{Proof} Since $\iota_*\circ\iota_!(a)=[c_0]\cdot a$ for $a\in
H_*(LM)$, the assertion follows from Corollary 3.3. 
\qed
\enddemo

Next, we show that similar statements hold for loop bracket products of
the form $\{[c_0], a\}$ for $a\in H_*(LM)$.

\proclaim{Corollary 3.5} Suppose $\chi(M)\not=0$ for a closed connected
oriented smooth $d$-manifold $M$, and let $a\in H_*(LM)$. 
\roster
\item If $|a|\not=d, d-1$, then the element $\{[c_0], a\}$ is either
$0$ or a torsion element of order a divisor of $\chi(M)$. 
\item Suppose further $M$ is simply connected. Then if $|a|\not=d-1$,
then the element $\{[c_0], a\}$ is either $0$ or a torsion element of
order a divisor of $\chi(M)$.
\endroster
\endproclaim
\demo{Proof} Since $\chi(M)\not=0$, $M$ is even dimensional. The
BV-identity multiplied by $\chi(M)$ gives 
$$
\Delta\bigl(\chi(M)[c_0]\cdot a\bigr)=\chi(M)\Delta([c_0])\cdot
a+\chi(M)[c_0]\cdot \Delta(a)+\chi(M)\{[c_0], a\}.
$$
If $|a|\not=d, d-1$, then by Corollary 3.3, we have $\chi(M)[c_0]\cdot
a=0$ and $\chi(M)[c_0]\cdot\Delta(a)=0$. Since $S^1$ action on $M$ is
trivial, we have $\Delta([c_0])=0$. Thus $\chi(M)\{[c_0],a\}=0$, and
the conclusion of (1) follows. 

For (2), when $|a|=d$, the element $\Delta(a)$ has degree $d+1$. By
Corollary 3.3, $\chi(M)[c_0]\cdot \Delta(a)=0$. If $M$ is simply
connected, $LM$ has a single component $L_0M$ and so $[c_0]\cdot a\in
H_0(LM)\cong \Bbb Z$ generated by $[c_0]$. Since $\Delta([c_0])=0$, we
have $\Delta([c_0]\cdot a)=0$. Hence $\chi(M)\{[c_0],a\}=0$, from
which the conclusion follows. 
\qed
\enddemo

When $|a|=d-1$, since $\{[c_0],a\}\in H_0(LM)$, this element is either
$0$ or torsion free. To see what happens when $M$ is not simply
connected, for each conjugacy class $[g]$ of $\pi_1(M)$ we choose a
loop $\gamma_g$ in $M$ belonging to $[g]$. When $|a|=d$, the element
$[c_0]\cdot a$ is a linear combination of classes $[\gamma_g]\in
H_0(LM)$. Since $\Delta([\gamma_g])\in H_1(L_{[g]}M)$ can be nonzero,
the simple connectivity assumption is needed in (2) of Corollary 3.5. 

\smallskip

{\bf Example 3.6}. We can verify Corollary 3.3 in actual examples. In
\cite{CJY}, the loop homology algebra for $LS^{2n}$ and $L\Bbb CP^n$
are computed. Their computation shows
$$
\align 
\Bbb H_*(LS^{2n})&\cong \Lambda(b)\otimes\Bbb
Z[a,v]/(a^2,ab,2av),\qquad
b\in\Bbb H_{-1}, a\in\Bbb H_{-2n}, v\in\Bbb H_{4n-2}, \\
\Bbb H_*(L\Bbb CP^n)&\cong \Lambda(w)\otimes\Bbb Z[c,u]/(c^{n+1},(n+1)c^nu,
wc^n), \quad w\in\Bbb H_{-1}, c\in\Bbb H_{-2}, u\in\Bbb H_{2n}.
\endalign
$$
For $\Bbb H_*(LS^{2n})$, we have $[c_0]=a$ and $\chi(S^{2n})=2$. By
the above computation, we can easily see that 
$\chi(S^{2n})[c_0]\cdot x=2a\cdot x=0$ for all
$x\in\Bbb H_*(LS^{2n})$ not in $\Bbb H_0$. For $\Bbb H_*(L\Bbb CP^n)$,
we have $[c_0]=c^n$ and $\chi(\Bbb CP^n)=n+1$. Again we can easily see
that the identity $\chi(\Bbb CP^n)[c_0]\cdot y=(n+1)c^n\cdot y=0$ 
for all $y$ not in $\Bbb H_0$. 

\smallskip

We discuss two final related results. The first one concerns an analogue of 
the BV identity for the loop coproduct. The BV identity can be understood by saying that the
failure of the commutativity of the following diagram is the loop
bracket:
$$
\CD
H_*(LM)\otimes H_*(LM) @>{\text{loop product}}>>  H_*(LM)  \\
@V{\Delta\otimes 1+1\otimes\Delta}VV  @V{\Delta}VV \\
H_*(LM)\otimes H_*(LM)  @>{\text{loop product}}>>   H_*(LM).
\endCD
$$
We ask a similar question for the loop coproduct. Does the following
diagram commute? If not, what is the measure of the failure of the
commutativity? 
$$
\CD
H_*(LM) @>{\Psi}>> H_*(LM\times LM) \\
@V{\Delta}VV  @V{\Delta\times 1+ 1\times\Delta}VV  \\
H_*(LM) @>{\Psi}>>  H_*(LM\times LM)
\endCD
$$
Unfortunately, things turn out to be rather trivial for the loop
coproduct. 

\proclaim{Proposition 3.7} For every $a\in H_*(LM)$, the identity
$(\Delta\times 1+1\times\Delta)\Psi(a)=0$ holds. 
\endproclaim
\demo{Proof} For $a\in H_*(LM)$, by (2) of Theorem 3.1, $\Psi(a)\in\Bbb Z[c_0]\otimes[c_0]\subset H_0(LM\times LM)$. Since $\Delta([c_0])=0$, the above identity holds. 
\qed
\enddemo

For the second result, recall that the loop product and the loop
coproduct satisfy Frobenius compatibility (Theorem 2.2). We ask a
similar question. What is the compatibility relation for the loop
bracket and the loop coproduct? The result turns out to be trivial when
one of the elements is from $H_*(M)$. 

\proclaim{Proposition 3.8} Let $M$ be as before with
$\chi(M)\ne0$. Suppose $a\in H_*(M)$. Then for any $b\in H_*(LM)$, we
have $\Psi(\{a,b\})=0$. 
\endproclaim
\demo{Proof} Let $\alpha\in H^*(M)$ be the cohomology class dual to $a$. 
Since $\Delta\alpha\cap b=(-1)^{|\alpha|}\{a,b\}$ (see \cite{T1}), using the
coderivation property of the cap product with respect to the loop
coproduct (2-5), 
$$
\align
\Psi(\{a,b\})&=(-1)^{|\alpha|+(|\alpha|-1)d}
(\Delta\alpha\times 1 + 1 \times \Delta\alpha)\cap\Psi(b)\\
&=(-1)^{|\alpha|+(|\alpha|-1)d}
\bigl[\chi(M)\bigl(\Delta\alpha\cap([c_0]\cdot b)\bigr)\otimes[c_0]
+\chi(M)[c_0]\otimes\bigl(\Delta\alpha\cap ([c_0]\cdot b)\bigr)\bigr].
\endalign
$$
Since the loop bracket behaves as a derivation in each variable, and $\{a,[c_0]\}=0$ for $a\in H_*(M)$, we have $\Delta\alpha\cap([c_0]\cdot
b)=(-1)^{|\alpha|}\{a,[c_0]\cdot b\}
=(-1)^{|\alpha|+(|\alpha|+1)d}[c_0]\cdot\{a,b\}$. The above identity then becomes 
$$
\Psi(\{a,b\})=\chi(M)([c_0]\cdot\{a,b\})\otimes[c_0] +
\chi(M)[c_0]\otimes([c_0]\cdot\{a,b\})=\Psi(\{a,b\})+\Psi(\{a,b\}),
$$
using (3-3). 
Hence $\Psi(\{a,b\})=0$. 
\qed
\enddemo

\enddocument